\documentclass{amsart}
\usepackage{amsmath,amssymb,verbatim,graphics,axodraw}

\numberwithin{equation}{section}

\newtheorem{prop}{Proposition}
\newtheorem{thm}[prop]{Theorem}

\newtheorem{lem}[prop]{Lemma}

\theoremstyle{definition}

\newtheorem{rem}[prop]{Remark}

\begin{document}

\title{Bosonic formula for level-restricted paths}

\author{Anne Schilling${}^*$}
\address{Instituut voor Theoretische Fysica\\
Universiteit van Amsterdam\\
Valckenierstraat 65\\
1018 XE Amsterdam\\
The Netherlands}
\email{schillin@wins.uva.nl}
\thanks{${}^*$ Supported by the 
``Stichting Fundamenteel Onderzoek der Materie''.}
\author{Mark Shimozono${}^\dag$}
\address{Department of Mathematics\\
Virginia Tech\\
Blacksburg, VA 24061-0123\\
U.S.A.}
\email{mshimo@math.vt.edu}
\thanks{${}^{\dag}$ Partially supported by NSF grant DMS-9800941.}

\date{December 1998}

\begin{abstract}
We prove a bosonic formula for the generating function of
level-restricted paths for the infinite families of
affine Kac-Moody algebras.  In affine type A this yields an
expression for the level-restricted generalized Kostka polynomials.
\end{abstract}

\maketitle

\newcommand{\la}{\lambda}
\newcommand{\La}{\Lambda}

\newcommand{\CC}{\mathbb{C}}
\newcommand{\Z}{\mathbb{Z}}
\newcommand{\ZNN}{\mathbb{Z}_{\ge0}}

\newcommand{\gggg}{\mathfrak{g}}
\newcommand{\gcl}{{\gggg}'}
\newcommand{\gfin}{\overline{\mathfrak{g}}}

\newcommand{\hhh}{\mathfrak{h}}
\newcommand{\ack}{\alpha^\vee}
\newcommand{\hck}{h^\vee}
\newcommand{\ac}{a^\vee}
\newcommand{\thck}{\theta^\vee}

\newcommand{\Paf}{P}
\newcommand{\Pcl}{P_{cl}}
\newcommand{\Pfin}{\overline{P}}
\newcommand{\Qfin}{\overline{Q}}

\newcommand{\Xfin}{\overline{X}}

\newcommand{\ba}[1]{\overline{#1}}

\newcommand{\boxf}[1]{
\raisebox{-0.15cm}{\scalebox{0.8}{\begin{picture}(17,20)(-2,-10)
\BBoxc(6.5,0)(13,13)
\Text(7,0)[]{\footnotesize$\mbox{#1}$}
\end{picture}}}}

\newcommand{\acl}{\alpha^{cl}}
\newcommand{\Lacl}{\La^{cl}}

\newcommand{\alb}{\overline{\alpha}}
\newcommand{\Lab}{\overline{\La}}

\newcommand{\af}{\mathrm{af}}
\newcommand{\cl}{\mathrm{cl}}

\newcommand{\rhoaf}{\rho}
\newcommand{\rhofin}{\overline{\rho}}

\newcommand{\Wfin}{\overline{W}}

\newcommand{\Aaf}{A_{\mathrm{af}}}
\newcommand{\Afin}{\overline{A}}

\newcommand{\Pifin}{\overline{\Pi}}
\newcommand{\Jfin}{\overline{J}}

\newcommand{\SSS}{\mathcal{S}}

\newcommand{\ch}{\mathrm{ch}}
\newcommand{\wt}{\mathrm{wt}}
\newcommand{\sign}{\mathrm{sign}}
\newcommand{\weight}{\mathrm{weight}}
\newcommand{\shape}{\mathrm{shape}}
\newcommand{\mult}{\mathrm{mult}}

\newcommand{\sln}{sl_n}
\newcommand{\slnhat}{\widehat{\sln}}

\newcommand{\inner}[2]{\langle #1,#2\rangle}
\newcommand{\Pu}{\mathcal{P}}
\newcommand{\VVV}{\mathbb{V}}
\newcommand{\Vfin}{V}
\newcommand{\BBB}{\mathbb{B}}

\newcommand{\KP}{K}
\newcommand{\Bf}{\mathcal{B}}

\section{Introduction}

Let $\gggg$ be an affine Kac-Moody algebra,
$V$ a $U_q(\gggg)^+$-submodule of a finite direct sum
$V'$ of irreducible integrable highest weight $U_q(\gggg)$-modules,
and $\Pi$ the limit of the Demazure operator for an element $w$
of the Weyl group as $\ell(w)\rightarrow\infty$.
The main theorem of this paper gives
sufficient conditions on $V$ so that the formula
\begin{equation} \label{first pi eq}
  \Pi\,\, \ch(V) = \ch(V')
\end{equation}
holds, where $\ch(V)$ is the character of $V$.
When $V$ is the one-dimen\-sional $U_q(\gggg)^+$-module
generated by the dominant integral weight $\La$ then \eqref{first pi eq}
is the Weyl-Kac character formula.  The above result
is well-known when $V$ is a union of Demazure modules
for any Kac-Moody algebra $\gggg$.

Let $\gcl$ be the derived subalgebra of $\gggg$.
Consider the $\gcl$-module $V$ given by a tensor product
of finite-dimensional $U_q(\gcl)$-modules that admit a crystal
of level at most $\ell$, 
with the one-dimensional subspace generated by a 
highest weight vector of an irreducible integrable
highest weight $U_q(\gcl)$-module of level $\ell$.  Such modules $V$
can be given the structure of a $U_q(\gggg)^+$-module
and as such, satisfy the above conditions.
Then a special case of \eqref{first pi eq} is a 
bosonic formula for the $q$-enumeration of level-restricted
inhomogeneous paths by the energy function.  In type
$A^{(1)}_{n-1}$ this formula was conjectured in \cite{FLOT},
stated there as a $q$-analogue of the Goodman-Wenzl straightening
algorithm for outer tensor products of irreducible modules
over the type $A$ Hecke algebra at a root of unity \cite{GW}.
In the isotypic component of the vacuum,
the bosonic formula coincides with half of the
bose-fermi conjecture in \cite[(9.2)]{SW}.  

The authors would like to thank Nantel Bergeron, Omar Foda, Masa\-ki Kashiwara,
Atsuo Kuniba, Masato Okado, Jean-Yves Thibon, and Ole Warnaar for helpful 
discussions.

\section{Notation}

Most of the following notation is taken from ref.~\cite{Kac}.
Let $X$ be a Dynkin diagram of affine type with vertices
indexed by the set $I=\{0,1,2,\dots,n\}$ as in \cite{Kac},
Cartan matrix $A=(a_{ij})_{i,j\in I}$, $\gggg=\gggg(A)$ the affine
Kac-Moody algebra, and $\hhh$ the Cartan subalgebra.
Let $\{\ack_i:i\in I\} \subset \hhh$ and $\{\alpha_j:j\in I\}\subset\hhh^*$
be the simple coroots and roots, which are linearly independent
subsets that satisfy $\inner{\ack_i}{\alpha_j}=a_{ij}$
for $i,j\in I$ where $\inner{\cdot}{\cdot}:\hhh\otimes\hhh^*\rightarrow\CC$
is the natural pairing.  Let $Q=\bigoplus_{i\in I} \Z \alpha_i$
be the root lattice.  Let the null root $\delta = \sum_{i\in I} a_i
\alpha_i$ be the unique element of the positive cone
$\bigoplus_{i\in I} \ZNN\, \alpha_i$ in $Q$,
that generates the one-dimensional lattice
$\{\beta\in Q | \inner{\ack_i}{\beta} = 0 \text{ for all $i\in I$ }\}$.
Let $K=\sum_{i\in I} \ac_i \ack_i\in \hhh$ be the canonical central element,
where the integers $\ac_i$ are the analogues of the integers $a_i$
for the dual algebra $\gggg^\vee$ defined by
the transpose ${}^tA$ of the Cartan matrix $A$.
Let $d\in \hhh$ (the degree derivation) be defined by the
conditions $\inner{d}{\alpha_i} = \delta_{i0}$ where $\delta_{ij}$ is
the Kronecker delta; $d$ is well-defined up to a summand proportional to
$K$.  Then $\{\ack_0,\dots,\ack_n,d\}$ is a basis of $\hhh$.
Let $\{\La_0,\dots,\La_n,\delta\}$ be the dual basis of $\hhh^*$;
the elements $\{\La_0,\dots,\La_n\}$ are called the fundamental weights.
The weight lattice is defined by
$P=\bigoplus_{i\in I} \Z \La_i \bigoplus \Z a_0^{-1} \delta$;
in the usual definition the scalar $a_0^{-1}$ is absent.
The weight lattice contains the root lattice since
$\alpha_j = \sum_{i\in I} a_{ij} \Lambda_i$ for $j\in I$.
Define $P^+ = \bigoplus_{i\in I} \ZNN\,\La_i \bigoplus \Z a_0^{-1} \delta$.
Say that a weight $\La\in P^+$ has level $\ell$ if
$\ell = \inner{K}{\La}$.

Consider the standard symmetric bilinear form $(\cdot|\cdot)$ on $\hhh^*$.
Since $\{\alpha_0,\dots,\alpha_n,\La_0\}$ is a basis of $\hhh^*$,
this form is uniquely defined by
$(\alpha_i|\alpha_j)=\ac_i a_i^{-1} a_{ij}$ for $i,j\in I$,
$(\alpha_i|\La_0)=\delta_{i0} a_0^{-1}$ for $i\in I$ and
$(\La_0|\La_0)=0$.  This form induces an isomorphism
$\nu:\hhh\rightarrow \hhh^*$ defined by
$\ac_i \nu(\ack_i) = a_i \alpha_i$ for $i\in I$
and $\nu(d)=a_0 \La_0$.  Also $\nu(K)=\delta$.

The Weyl group $W$ is the subgroup of $GL(\hhh^*)$
generated by the simple reflections $r_i$ ($i\in I$) defined by
\begin{equation*}
  r_i (\beta) = \beta - \inner{\ack_i}{\beta} \alpha_i.
\end{equation*}
The form $(\cdot|\cdot)$ is $W$-invariant.
Suppose $\alpha\in Q$ is a real root, that is, the $\alpha$-weight
space of $\gggg$ is nonzero and there is a simple root $\alpha_i$ and a
Weyl group element $w\in W$ such that $\alpha=w(\alpha_i)$.
Define $\alpha^\vee \in \hhh$ by $w(\ack_i)$.  This is independent
of the expression $\alpha = w(\alpha_i)$.  Define
$r_\alpha \in W$ by 
\begin{equation*}
  r_\alpha(\beta) = \beta - \inner{\alpha^\vee}{\beta} \alpha\qquad
  \text{for $\beta\in\hhh^*$.}
\end{equation*}

Let $\gcl$ be the derived algebra of $\gggg$, obtained by
``omitting" the degree derivation $d$.  Its weight lattice
is $\Pcl \cong P/\Z a_0^{-1}\delta$.  Denote the canonical
projection $P \rightarrow \Pcl$ by $\cl$.
Write $\acl_i=\cl(\alpha_i)$ and
$\Lacl_i=\cl(\La_i)$ for $i\in I$.
The elements $\{\acl_i \mid i\in I\}$ are linearly dependent.
Write $\af:\Pcl\rightarrow P$ for the section of $\cl$
given by $\af(\Lacl_i)=\La_i$ for all $i\in I$.
Write $\Pcl^+ = \bigoplus_{i\in I} \ZNN\, \Lacl_i$.
Define the level of $\mu\in \Pcl^+$
to be $\inner{K}{\af(\mu)}$.

Consider the Dynkin diagram $\Xfin$ obtained by removing the vertex $0$
from the diagram $X$, with corresponding Cartan matrix
$\Afin$ indexed by the set $J=I-\{0\}$, and let
$\gfin=\gggg(\Afin)$ be the simple Lie algebra.
One has the inclusions $\gfin \subset \gcl \subset \gggg$.
Let $\{\alb_i: i\in J\}$ be the simple roots,
$\{\Lab_i: i\in J\}$ the fundamental weights, and
$\Qfin=\bigoplus_{i\in J} \Z \alb_i$ the root lattice for $\gfin$.
The weight lattice of $\gfin$ is $\Pfin=\bigoplus_{i\in J} \Z\Lab_i$
and $\Pfin \cong \Pcl/\Z\La_0$.  The image of $\La\in P$ into
$\Pfin$ is denoted by $\overline{\La}$.
We shall use the section of the natural projection
$\Pcl\rightarrow \Pfin$ given by the map $\Pfin\rightarrow\Pcl$
that sends $\Lab_i \mapsto \Lacl_i-\Lacl_0$ for $i\in J$.
By abuse of notation, for $\La\in P$, $\overline{\La}$
shall also denote the image of the element $\overline{\La}$
under the lifting map $\Pfin \rightarrow P$ specified above.

Let $\Pfin^+=\bigoplus_{i\in J} \ZNN \,\Lab_i$.
For $\la \in \Pfin^+$, denote by $\Vfin(\la)$ the
irreducible integrable highest weight $U_q(\gfin)$-module
of highest weight $\la$.

Let $\theta = \delta - a_0 \alpha_0 = \sum_{i\in J} a_i \alpha_i \in \Qfin$.
One has the formulas $(\theta|\theta)=2 a_0$, $\theta=a_0 \nu(\thck)$, and
$\ack_0 = K - a_0 \thck$.  Observe that
\begin{equation*}
  \cl(\alpha_0)=-a_0^{-1}\sum_{i\in J} a_i \acl_i=-\cl(\nu(\thck)).
\end{equation*}

For $\La\in P^+$ let $\VVV(\La)$ be the
irreducible integral highest weight module of highest weight $\La$
over the quantized universal enveloping algebra $U_q(\gggg)$,
$\BBB(\La)$ the crystal base of $\VVV(\La)$,
and $u_\La\in \BBB(\La)$ the highest weight vector.

By restriction from $U_q(\gggg)$ to $U_q(\gcl)$,
the module $\VVV(\La)$ is an irreducible integral highest
weight module for $U_q(\gcl)$ of highest weight $\cl(\La)$,
with crystal $\BBB(\La)$ that is $\Pcl$-weighted by composing the
weight function $\BBB(\La)\rightarrow P$ with the projection $\cl$.
Conversely, any integrable irreducible highest weight $U_q(\gcl)$-module
can be obtained this way.

\section{Short review of affine crystal theory}

\subsection{Crystals}

A $P$-weighted $I$-crystal $B$ is a colored graph with
vertices indexed by $b\in B$, directed edges colored by $i\in I$,
and a weight function $\wt:B\rightarrow P$, satisfying the 
axioms below.  First some notation is required.  Denote an
edge from $b$ to $b'$ colored $i$, by $b' = f_i(b)$
or equivalently $b=e_i(b')$.  Write $\phi_i(b)$
(resp. $\epsilon_i(b)$) for the maximum index $m$ such
that $f_i^m(b)$ (resp. $e_i^m(b)$) is defined.
\begin{enumerate}
\item If $b'=f_i(b)$ then $\wt(b')=\wt(b)-\alpha_i$.
\item $\phi_i(b)-\epsilon_i(b) = \inner{\ack_i}{\wt(b)}$.
\end{enumerate}
An element $u\in B$ is a highest weight vector if
$e_i(u)$ is undefined for all $i\in I$.
The $i$-string of $b\in B$ consists of all elements
$e_i^m(b)$ ($0\le m\le \epsilon_i(b)$) and
$f_i^m(b)$ ($0\le m\le \phi_i(b)$).
The nondominant part of the $i$-string is comprised of
all elements which admit $e_i$.

We also define the crystal reflection operator $s_i:B\rightarrow B$ by
\begin{equation*}
  s_i(b) = \begin{cases}
  f_i^{\phi_i(b)-\epsilon_i(b)}(b) & \text{if $\phi_i(b)>\epsilon_i(b)$} \\
  b & \text{if $\phi_i(b)=\epsilon_i(b)$} \\
  e_i^{\epsilon_i(b)-\phi_i(b)}(b) & \text{if $\phi_i(b)<\epsilon_i(b)$}.
	\end{cases}
\end{equation*}
It is obvious that $s_i$ is an involution.  Observe that
\begin{equation} \label{refl wt}
  \wt(s_i(b))=r_i\wt(b)=\wt(b) - \inner{\ack_i}{\wt(b)}\alpha_i.
\end{equation}

Define the notation $\phi(b)=\sum_{i\in I} \phi_i(b) \La_i$
and $\epsilon(b)=\sum_{i\in I} \epsilon_i(b) \La_i$.

If a $U_q(\gggg)$-module (resp. $U_q(\gcl)$-module, resp.
$U_q(\gfin)$-module) has a crystal base then the latter is
naturally a $P$-weighted (resp. $\Pcl$-weighted, resp. $\Pfin$-weighted)
$I$-crystal (resp. $I$-crystal, resp. $J$-crystal).

\subsection{Tensor products}

Given crystals $B_1$ and $B_2$, contrary to the
literature (but consistent with the Robinson-Schensted-Knuth
correspondence in type A), define the following crystal
structure on the tensor product $B_2 \otimes B_1$.
The elements are denoted $b_2 \otimes b_1$ for $b_i\in B_i$
($i\in\{1,2\}$) and one defines
\begin{equation*}
\begin{split}
  \phi_i(b_2 \otimes b_1) &= \phi_i(b_2) + \max(0,\phi_i(b_1)-\epsilon_i(b_2)) \\
  \epsilon_i(b_2 \otimes b_1) &= \epsilon_i(b_1) + \max(0,-\phi_i(b_1)+\epsilon_i(b_2)).
\end{split}
\end{equation*}
When $\phi_i(b_2 \otimes b_1) > 0$ (resp. $\epsilon_i(b_2 \otimes b_1) > 0$)
one defines
\begin{align*}
  f_i(b_2 \otimes b_1) &= \begin{cases}
	b_2 \otimes f_i(b_1) & \text{if $\phi_i(b_1)>\epsilon_i(b_2)$} \\
	f_i(b_2) \otimes b_1 & \text{if $\phi_i(b_1)\le\epsilon_i(b_2)$}
	\end{cases} \\
\intertext{and respectively}
  e_i(b_2 \otimes b_1) &= \begin{cases}
	b_2 \otimes e_i(b_1) & \text{if $\phi_i(b_1)\ge\epsilon_i(b_2)$}\\
	e_i(b_2) \otimes b_1 & \text{if $\phi_i(b_1)<\epsilon_i(b_2)$}.
  \end{cases}
\end{align*}

An element of a tensor product of crystals is called a path.

\subsection{Energy function}

The definitions here follow \cite{NY}.
Suppose that $B_1$ and $B_2$ are crystals of finite-dimensional
$U_q(\gcl)$-modules such that $B_2 \otimes B_1$ is connected.
Then there is an isomorphism of $\Pcl$-weighted $I$-crystals
$B_2 \otimes B_1 \cong B_1 \otimes B_2$.  This is called
the local isomorphism.  Let the image of
$b_2\otimes b_1\in B_2\otimes B_1$ under this isomorphism
be denoted $b_1' \otimes b_2'$.  Then there is a unique
(up to a global additive constant) map
$H:B_2 \otimes B_1 \rightarrow \Z$ such that
\begin{equation*}
  H(e_i(b_2 \otimes b_1)) = H(b_2\otimes b_1) +
  \begin{cases}
    -1 & \text{if $i=0$, $e_0(b_2\otimes b_1)=e_0(b_2) \otimes b_1$} \\
       & \text{and $e_0(b_1'\otimes b_2')=e_0(b_1') \otimes b_2',$} \\
    1 & \text{if $i=0$, $e_0(b_2\otimes b_1)=b_2 \otimes e_0(b_1)$} \\
    & \text{and $e_0(b_1'\otimes b_2')=b_1' \otimes e_0(b_2'),$} \\
    0 & \text{otherwise.}
  \end{cases}
\end{equation*}
This map is called the local energy function.

Consider $B = B_L \otimes \dots \otimes B_1$ with
$B_j$ the crystal of a finite-dimensional $U_q(\gcl)$-module
for $1\le j\le L$.  Assume that for all $1\le i<j\le L$,
$B_j \otimes B_i$ is a connected $\Pcl$-weighted $I$-crystal.
Given $b=b_L \otimes \dots \otimes b_1\in B$,
denote by $b_j^{(i+1)}$ the $(i+1)$-th tensor factor in the
image of $b$ under the composition of local isomorphisms that switch
$B_j$ with $B_k$ as $k$ goes from $j-1$ down to $i+1$.
Then define the energy function
\begin{equation} \label{energy}
  E_B(b) = \sum_{1\le i<j\le L} H_{j,i}(b_j^{(i+1)}\otimes b_i)
\end{equation}
where $H_{j,i}:B_j \otimes B_i\rightarrow\Z$ is the local
energy function.  It satisfies the following property.

\begin{lem} \label{energy e} \cite[Prop. 1.1]{HKKOTY}
Suppose $i\in I$, $b\in B$ and $e_i(b)$ is defined.
If $i\not=0$ then $E_B(e_i(b))=E_B(b)$.  If $i=0$ and
$b$ has the property that for any of its images
$b'=b'_L \otimes \dots \otimes b'_1$ under a composition
of local isomorphisms, $e_0(b')=b'_L \otimes
\dots \otimes e_0(b'_k) \otimes \dots \otimes b'_1$ with
$k\not=1$, then $E_B(e_0(b))=E_B(b)-1$.
\end{lem}

\subsection{Classically restricted paths}

Say that $b\in B:=B_L\otimes\dots\otimes B_1$ is classically restricted if
$b$ is a $\gfin$-highest weight vector, that is, $e_i(b)$ is undefined
for all $i\in J$.  For $\la\in\Pfin^+$ denote by
$\Pu(B,\la)$ the set of classically restricted $b\in B$
of weight $\la$.  Define the polynomial
\begin{equation}
  \KP(B,\la)(q) = \sum_{b\in \Pu(B,\la)} q^{E_B(b)}
\end{equation}
where $E_B$ is the energy function on $B$.  For $\gggg$ of type $A^{(1)}_{n-1}$
$\KP(B,\la)(q)$ is the generalized Kostka polynomial \cite{S1,S2,SW}.

\subsection{Almost perfect crystals}

Let $B$ be the crystal of a finite-dimen\-sional $U_q(\gcl)$-module.
Say that $B$ is almost perfect of level $\ell$ \cite{Ok} if
it satisfies the following weakening of the definition
of a perfect crystal \cite[Def. 4.6.1]{KMN}:
\begin{enumerate}
\item $B \otimes B$ is connected.
\item There is a $\La'\in\Pcl$ such that there is a unique
$b'\in B$ such that $\wt(b')=\La'$ and for every $b\in B$,
$\wt(b) \in \La' - \bigoplus_{i\in J} \ZNN\alpha_i$.
\item For every $b\in B$, $\inner{K}{\epsilon(b)} \ge \ell$.
\item For every $\La\in\Pcl^+$ of level $\ell$, there is a $b,b'\in B$
such that $\epsilon(b)=\phi(b')=\La$.
\end{enumerate}
$B$ is said to be perfect if the elements $b$ and $b'$ in
item 4 are unique.

\subsection{Level restricted paths}
{}From now on, fix a positive integer $\ell$ (the level).

For $1\le j\le L$ let $B_j$ be the crystal of a 
finite-dimensional $U_q(\gcl)$-module, that is almost perfect
of level at most $\ell$.

Let $B=B_L \otimes \dots \otimes B_1$,
$\La,\La'\in \Pcl^+$ weights of level $\ell$, and
$\Pu(B,\La,\La')$ the set of paths
$b=b_L \otimes\dots\otimes b_1\in B$
such that $b \otimes u_\La \in B \otimes \BBB(\La)$
is a highest weight vector of weight $\La'$.

In the special case that $\La=\ell \La_0$,
the elements of $\Pu(B,\La,\La')$ are called the
level-$\ell$ restricted paths of weight $\La'$.

\begin{thm} \label{level tensor} \cite{KMN} \cite[Appendix A]{KMPY}.
Let $\gggg$ be an affine Kac-Moody algebra in one of the
infinite families. Let $B$ be the tensor product
of crystals of finite-dimensional $U_q(\gcl)$-modules, that
are almost perfect of level at most $\ell$, and
$\La\in\Pcl^+$ a weight of level $\ell$.  Then
there is an isomorphism of $\Pcl$-weighted $I$-crystals 
\begin{equation}\label{tensor iso}
  B \otimes \BBB(\La) \cong
  \bigoplus_{\La'\in \Pcl^+} 
  \bigoplus_{b\in \Pu(B,\La,\La')} \BBB(\La')
\end{equation}
where $\La'$ is of level $\ell$.
\end{thm}

This isomorphism of $\Pcl$-weighted crystals can be
lifted to one of $P$-weighted crystals by specifying an integer multiple
of $a_0^{-1}\delta$ for each highest weight vector in $B\otimes \BBB(\La)$.
However for our purposes this should be done in a way that extends
the definition of the energy function for $B$.  To this end,
choose a perfect crystal $B_0$ of level $\ell$, and 
assume that for all $0\le i<j\le L$, $B_j\otimes B_i$ is
connected.
Let $b_0\in B_0$ be the unique element such that $\phi(b_0)=\La$.
Define the energy function
$E:B\rightarrow \Z$ by $E(b)=E_{B,B_0}(b \otimes b_0)$
where $E_{B,B_0}:B\otimes B_0\rightarrow \Z$ is the energy function
defined in \eqref{energy}.
For $b\in \Pu(B,\La,\La')$, define an affine weight function
$\wt(b\otimes u_\La)=\af(\La')-E(b)a_0^{-1}\delta$.
This defines the $P$-weight of every highest weight vector in
$B\otimes \BBB(\La)$ and hence a $P$-weight function for all of
$B\otimes \BBB(\La)$.

Then one has the following $P$-weighted analogue of
\eqref{tensor iso}:
\begin{equation}\label{wt tensor iso}
  B \otimes \BBB(\af(\La)) \cong
  \bigoplus_{\La'\in \Pcl^+}
  \bigoplus_{b\in \Pu(B,\La,\La')} \BBB(\wt(b\otimes u_\La))
\end{equation}
where $\La'$ is of level $\ell$.  This decomposition
can be described by the polynomial
\begin{equation}\label{K def}
  \KP(B,\La,\La',B_0)(q) = \sum_{b\in \Pu(B,\La,\La')} q^{E(b)}.
\end{equation}
Our goal is to prove a formula for the polynomial $\KP(B,\La,\La',B_0)(q)$.

\section{General bosonic formula}
\label{sec general bose}

Let $J$ be the antisymmetrizer
\begin{equation*}
  J = \sum_{w\in W} (-1)^w w.
\end{equation*}
Write
\begin{equation*}
  R = \prod_{\alpha\in\Delta_+} (1-\exp(-\alpha))^{\mult(\alpha)}
\end{equation*}
where $\Delta_+$ is the set of positive roots of $\gggg$
and $\mult(\alpha)$ is the dimension of the $\alpha$-weight
space in $\gggg$.

Let $\rho\in P^+$ be the unique weight defined by 
$\inner{\ack_i}{\rho}=1$ for all $i\in I$
and $\inner{d}{\rho}=0$.  It satisfies
$\inner{\thck}{\rho}=a_0^{-1} \inner{K-\ack_0}{\rho}=a_0^{-1}(\hck-1)$
where $\hck=\sum_{i\in I}\ac_i$ is the dual Coxeter number.
Define the operator
\begin{equation*}
  \Pi(p) = R^{-1} e^{-\rho} J(e^\rho p).
\end{equation*}
where $R^{-1}$ makes sense by expanding the reciprocals of the
factors of $R$ in geometric series.
The computation is defined in a suitable completion of $\Z[P]$.
One has $\Pi\,(e^\La)=\ch\, \VVV(\La)$ for all $\La\in P^+$,
which is the Weyl-Kac character formula \cite[Theorem 10.4]{Kac}.

\begin{thm} \label{gen boson} Let $\gggg$
be an affine Kac-Moody algebra
in one of the infinite families,
$B'$ the crystal of a finite direct sum of
irreducible integrable highest weight $U_q(\gggg)$-modules
and $B \subset B'$ a subset such that:
\begin{enumerate}
\item $B$ is closed under $e_i$ for all $i\in I$.
\item $B'$ is generated by $B$.
\item For all $b\in B$ and $i\in I$, if
$\epsilon_i(b)>0$ then the $i$-string of $b$ in $B'$
is contained in $B$.
\end{enumerate}
Then
\begin{equation} \label{pi eq}
  \Pi\,\, \ch(B) = \ch(B').
\end{equation}
\end{thm}
\begin{proof} Without loss of generality it may be assumed that
$B'=\BBB(\La)$ for some $\La\in P^+$.  Multiplying both sides of
\eqref{pi eq} by $R\, e^\rho$, one obtains
\begin{equation*}
  \sum_{(w,b)\in W\times B} (-1)^w w(e^{\wt(b)+\rho}) =
  \sum_{w\in W} (-1)^w w(e^{\La+\rho}).
\end{equation*}
Observe that both sides are $W$-alternating.  The $W$-alternants
have a basis given by $J(\La+\rho)$ where $\La\in P^+$.
Taking the coefficient of $e^{\La+\rho}$ on both sides, 
\begin{equation} \label{cancel eq}
  \sum_{(w,b)\in \SSS} (-1)^w = 1
\end{equation}
where $\SSS$ is the set of pairs $(w,b)\in W \times B$ such that
\begin{equation} \label{wt cond}
  \wt(b) = w^{-1}(\La+\rho)-\rho.
\end{equation}
Observe that if $(w,b)\in\SSS$ is such that $b$ is a highest weight vector,
then $w=1$ and $b=u_\La$, for both of the regular dominant
weights $\wt(b)+\rho$ and $\La+\rho$ are in the same $W$-orbit and hence
must be equal.  Conditions 1 and 2 ensure that $u_\La\in B$.
Let $\SSS'=\SSS-\{(1,u_\La)\}$.  It is enough to show that
there is an involution $\Phi:\SSS'\rightarrow \SSS'$ with no fixed points,
such that if $\Phi(w,b)=(w',b')$ then $w$ and $w'$ have opposite sign.
In this case $\Phi$ is said to be sign-reversing.
Let $\SSS_i$ be the set of pairs $(w,b)\in\SSS'$ such that
$\epsilon_i(b)>0$.  Define the map $\Phi_i:\SSS_i\rightarrow\SSS_i$ by
$\Phi_i(w,b)=(w r_i, s_i e_i(b))$.  Note that
$s_i e_i(b)\in B$ by condition 3.
The condition \eqref{wt cond} for $\Phi_i(w,b)$ is
\begin{equation*}
\begin{split}
  (w r_i)^{-1}(\La+\rho)-\rho &=
  r_i w^{-1}(\La+\rho) - r_i \rho + r_i \rho - \rho \\
  &= r_i(w^{-1}(\La+\rho)-\rho) -\inner{\ack_i}{\rho}\alpha_i \\
  &= r_i(\wt(b)) - \alpha_i = \wt(f_i s_i(b)) = \wt(s_i e_i(b)).
\end{split}
\end{equation*}
Since $s_i e_i(b)=f_i s_i(b)$, $\epsilon_i(s_i e_i(b))>0$,
so that $(w r_i,s_i e_i(b))\in \SSS_i$.  This shows that
$\Phi_i$ is well-defined.  It follows directly from the
definitions that $\Phi_i$ is a sign-reversing involution.

Since $\SSS' = \bigcup_{i\in I} \SSS_i$ it suffices to define
a global involutive choice of the canceling root direction for each pair
$(w,b)\in \SSS'$, that is, a function
$v:\SSS'\rightarrow I$ such that if $v(w,b)=i$ then
\begin{enumerate}
\item[(V1)] $(w,b)\in \SSS_i$.
\item[(V2)] $v(w r_i,s_i e_i(b))=i$.
\end{enumerate}

Let $\La = \La_{i_1}+\dots+\La_{i_\ell}$ be an expression of $\La$
as a sum of fundamental weights.  By \cite[Lemma 8.3.1]{Ka},
$\BBB(\La)$ is isomorphic to the full subcrystal of
$\BBB(\La_{i_\ell}) \otimes \dots \otimes \BBB(\La_{i_1})$
generated by $u_{\La_{i_\ell}} \otimes \dots \otimes u_{\La_{i_1}}$.

Given $(w,b)\in \SSS'$, let $b_\ell \otimes \dots \otimes b_1$
be the image of $b$ in the above tensor product of crystals
of modules of fundamental highest weight.
Let $r$ be minimal such that
$b_r \otimes b_{r-1} \otimes \dots \otimes b_1$ is not
a highest weight vector.  Then $b_{r-1} \otimes \dots \otimes b_1$
is a highest weight vector, say of weight $\La'$.

Let $\Bf$ be a perfect crystal of 
the same level as $\La_{i_r}$.
Given any $L>0$, the theory of perfect crystals 
\cite[Section 4.5]{KMN} gives an isomorphism of $P$-weighted crystals
\begin{equation*}
  \BBB(\La_{i_r}) \cong {\Bf}^{\otimes L} \otimes \BBB(\La_j)
\end{equation*}
where $j$ is determined by $i_r$ and $L$ and
$\Bf^{\otimes L}$ is $P$-weighted using the energy function.

Let $b_r\in \BBB(\La_{i_r})$ have image
$p_{-1} \otimes \dots \otimes p_{-L} \otimes u'$ where $u'
\in \BBB(\La_j)$.  Assume that $L$ is large enough so that $u'=u_{\La_j}$.
If one takes the image of $b_r$ in such a tensor product for $L'>L$
the tensor factors $p_{-1}$ through $p_{-L}$ do not change.

Let $k$ be minimal such that
$p_k \otimes \dots \otimes p_{-L} \otimes u_{\La_j} \otimes u_{\La'}$
is not a highest weight vector.  Observe that $k$ is independent of $L$
as long as $L$ is big enough.
Then $p_{k-1} \otimes \dots \otimes p_{-L} \otimes u_{\La_j} \otimes u_{\La'}$
is a highest weight vector, say of weight $\La''$.

So $p_k \in \Bf$ is such that
$\epsilon_i(p_k)>\inner{\ack_i}{\La''}$ for some $i\in I$;
let $I'$ be the set of such $i\in I$.

Fix an $i\in I'$.  Consider the same constructions for $b'=s_i e_i(b)$.
Let $b'_\ell \otimes \dots \otimes b'_1$ be the image of
$b'$ in the above tensor product of irreducible crystals of fundamental
highest weights.  Then $b'_{r-1} \otimes \dots\otimes b'_1 =
b_{r-1} \otimes \dots \otimes b_1$ and $b'_r \otimes \dots \otimes b'_1$
is not a highest weight vector; in particular it admits $e_i$.
Take L large enough so that the image of $b'_r$ in
${\Bf}^{\otimes L} \otimes \BBB(\La_j)$ also has the form
$p'_{-1} \otimes \dots \otimes p'_{-L} \otimes u_{\La_j}$.
Then $p_{k-1} \otimes \dots \otimes p_{-L} =
p'_{k-1} \otimes \dots \otimes p'_{-L}$ and
$p'_k \otimes \dots \otimes p'_{-L} \otimes u_{\La_j} \otimes u_{\La'}$
admits $e_i$.

The level of the fundamental weight $\La_i$ is $\ac_i$.
For the affine algebras $A_n^{(1)}$ and $C_n^{(1)}$,
$\ac_i=1$ for all $i\in I$. For all others $1\le \ac_i\le 2$.
The theorem now follows from Lemma \ref{lev one} below, applied with
the dominant integral weight $\La''$.
\end{proof}

\begin{lem} \label{lev one} 
For the affine Kac-Moody algebra $\gggg$
in one of the infinite families,
there exist perfect crystals $B$ of level one and two (the latter 
case only for $\gggg\neq A_n^{(1)}, C_n^{(1)}$) having the following 
property.
Let $\La$ be a dominant integral weight of positive level,
$S$ the set of elements $b_1\in B$ such that
$b_1 \otimes u_\La$ is not a highest weight vector.
Then there is a map $v:S\rightarrow I$ such that
if $v(b_1)=i$ then
\begin{enumerate}
\item $\epsilon_i(b_1\otimes u_\La)>0$.
\item For any $b_L,\dots,b_2\in B$, writing
$b'_L \otimes \dots \otimes b'_2 \otimes b'_1 \otimes u_\La
= s_i e_i(b_L \otimes \dots \otimes b_2 \otimes b_1 \otimes u_\La)$,
one has $b'_1 \in S$ and $v(b'_1)=i$.
\end{enumerate}
\end{lem}
\begin{proof} 
For the involutive property 2, it is sufficient that $v$ is constant
on the nondominant part of every string. Hence one only needs to
consider 
\begin{equation}\label{cond}
\text{$b_1$ that are on the nondominant part of at least two strings
of length $\ge 2$}.
\end{equation}

Perfect crystals of level one for $A_n^{(1)}$ ($n\ge 1$), 
$B_n^{(1)}$ ($n\ge 3$), $D_n^{(1)}$ ($n\ge 4$), $A_{2n}^{(2)}$ ($n\ge 1$),
$A_{2n-1}^{(2)}$ ($n\ge 3$) and $D_{n+1}^{(2)}$ ($n\ge 2$)
are listed in Table~\ref{table crystals} (see \cite[Section 6]{KMN}).
Note that there are no elements satisfying \eqref{cond}.
This guarantees the existence of the map $v$ with the desired properties.
\begin{table}
\begin{tabular}{|c|l|}
\hline
%
$A_n^{(1)}$ &
\raisebox{-0.7cm}{\scalebox{0.8}{
\begin{picture}(250,62)(-10,-12)
\BText(0,0){1}
\LongArrow(10,0)(40,0)
\BText(50,0){2}
\LongArrow(60,0)(90,0)
\BText(100,0){3}
\LongArrow(110,0)(140,0)
\Text(160,0)[]{$\cdots$}
\LongArrow(175,0)(205,0)
\BText(220,0){n+1}
\LongArrowArc(110,-181)(216,62,118)
\PText(25,2)(0)[b]{1}
\PText(75,2)(0)[b]{2}
\PText(125,2)(0)[b]{3}
\PText(190,2)(0)[b]{n}
\PText(110,38)(0)[b]{0}
\end{picture}
}}
\\ \hline
%
$B_n^{(1)}$ &
\raisebox{-1.3cm}{\scalebox{0.8}{
\begin{picture}(350,100)(-10,-50)
\BText(0,0){1}
\LongArrow(10,0)(30,0)
\BText(40,0){2}
\LongArrow(50,0)(70,0)
\Text(85,0)[]{$\cdots$}
\LongArrow(95,0)(115,0)
\BText(125,0){n}
\LongArrow(135,0)(155,0)
\BText(165,0){0}
\LongArrow(175,0)(195,0)
\BBoxc(205,0)(13,13)
\Text(205,0)[]{\footnotesize$\overline{\mbox{n}}$}
\LongArrow(215,0)(235,0)
\Text(250,0)[]{$\cdots$}
\LongArrow(260,0)(280,0)
\BBoxc(290,0)(13,13)
\Text(290,0)[]{\footnotesize$\overline{\mbox{2}}$}
\LongArrow(300,0)(320,0)
\BBoxc(330,0)(13,13)
\Text(330,0)[]{\footnotesize$\overline{\mbox{1}}$}
\LongArrowArc(185,-330)(365,68,112)
\LongArrowArcn(145,330)(365,-68,-112)
\PText(20,2)(0)[b]{1}
\PText(60,2)(0)[b]{2}
\PText(105,2)(0)[b]{n-1}
\PText(145,2)(0)[b]{n}
\PText(185,2)(0)[b]{n}
\PText(225,2)(0)[b]{n-1}
\PText(270,2)(0)[b]{2}
\PText(310,2)(0)[b]{1}
\PText(185,38)(0)[b]{0}
\PText(145,-35)(0)[t]{0}
\end{picture}
}}
\\ \hline
%
$D_n^{(1)}$ &
\raisebox{-1.3cm}{\scalebox{0.8}{
\begin{picture}(365,100)(-10,-50)
\BText(0,0){1}
\LongArrow(10,0)(30,0)
\BText(40,0){2}
\LongArrow(50,0)(70,0)
\Text(85,0)[]{$\cdots$}
\LongArrow(95,0)(115,0)
\BText(130,0){n-1}
\LongArrow(143,2)(160,14)
\LongArrow(143,-2)(160,-14)
\BText(170,15){n}
\BBoxc(170,-15)(13,13)
\Text(170,-15)[]{\footnotesize$\overline{\mbox{n}}$}
\LongArrow(180,14)(197,2)
\LongArrow(180,-14)(197,-2)
\BBoxc(215,0)(25,13)
\Text(215,0)[]{\footnotesize$\overline{\mbox{n-1}}$}
\LongArrow(230,0)(250,0)
\Text(265,0)[]{$\cdots$}
\LongArrow(275,0)(295,0)
\BBoxc(305,0)(13,13)
\Text(305,0)[]{\footnotesize$\overline{\mbox{2}}$}
\LongArrow(315,0)(335,0)
\BBoxc(345,0)(13,13)
\Text(345,0)[]{\footnotesize$\overline{\mbox{1}}$}
\LongArrowArc(192.5,-367)(402,69,111)
\LongArrowArcn(152.5,367)(402,-69,-111)
\PText(20,2)(0)[b]{1}
\PText(60,2)(0)[b]{2}
\PText(105,2)(0)[b]{n-2}
\PText(152,13)(0)[br]{n-1}
\PText(152,-9)(0)[tr]{n}
\PText(188,13)(0)[bl]{n}
\PText(188,-9)(0)[tl]{n-1}
\PText(240,2)(0)[b]{n-2}
\PText(285,2)(0)[b]{2}
\PText(325,2)(0)[b]{1}
\PText(192.5,38)(0)[b]{0}
\PText(152.5,-35)(0)[t]{0}
\end{picture}
}}
\\ \hline
%
$A_{2n}^{(2)}$ &
\raisebox{-0.7cm}{\scalebox{0.8}{
\begin{picture}(350,62)(-10,-12)
\BText(0,0){1}
\LongArrow(10,0)(30,0)
\BText(40,0){2}
\LongArrow(50,0)(70,0)
\Text(85,0)[]{$\cdots$}
\LongArrow(95,0)(115,0)
\BText(125,0){n}
\LongArrow(135,0)(155,0)
\BBoxc(165,0)(13,13)
\Text(165,0)[]{\footnotesize$\overline{\mbox{n}}$}
\LongArrow(175,0)(195,0)
\Text(210,0)[]{$\cdots$}
\LongArrow(220,0)(240,0)
\BBoxc(250,0)(13,13)
\Text(250,0)[]{\footnotesize$\overline{\mbox{2}}$}
\LongArrow(260,0)(280,0)
\BBoxc(290,0)(13,13)
\Text(290,0)[]{\footnotesize$\overline{\mbox{1}}$}
\Text(145,35)[]{\footnotesize$\emptyset$}
\LongArrowArc(145,-317)(352,67.5,89)
\LongArrowArc(145,-317)(352,91,112.5)
\PText(20,2)(0)[b]{1}
\PText(60,2)(0)[b]{2}
\PText(105,2)(0)[b]{n-1}
\PText(145,2)(0)[b]{n}
\PText(185,2)(0)[b]{n-1}
\PText(230,2)(0)[b]{2}
\PText(270,2)(0)[b]{1}
\PText(220,30)(0)[b]{0}
\PText(70,30)(0)[b]{0}
\end{picture}
}}
\\ \hline
%
$A_{2n-1}^{(2)}$ &
\raisebox{-1.3cm}{\scalebox{0.8}{
\begin{picture}(310,100)(-10,-50)
\BText(0,0){1}
\LongArrow(10,0)(30,0)
\BText(40,0){2}
\LongArrow(50,0)(70,0)
\Text(85,0)[]{$\cdots$}
\LongArrow(95,0)(115,0)
\BText(125,0){n}
\LongArrow(135,0)(155,0)
\BBoxc(165,0)(13,13)
\Text(165,0)[]{\footnotesize$\overline{\mbox{n}}$}
\LongArrow(175,0)(195,0)
\Text(210,0)[]{$\cdots$}
\LongArrow(220,0)(240,0)
\BBoxc(250,0)(13,13)
\Text(250,0)[]{\footnotesize$\overline{\mbox{2}}$}
\LongArrow(260,0)(280,0)
\BBoxc(290,0)(13,13)
\Text(290,0)[]{\footnotesize$\overline{\mbox{1}}$}
\LongArrowArc(165,-240)(275,65,115)
\LongArrowArcn(125,240)(275,-65,-115)
\PText(20,2)(0)[b]{1}
\PText(60,2)(0)[b]{2}
\PText(105,2)(0)[b]{n-1}
\PText(145,2)(0)[b]{n}
\PText(185,2)(0)[b]{n-1}
\PText(230,2)(0)[b]{2}
\PText(270,2)(0)[b]{1}
\PText(165,38)(0)[b]{0}
\PText(125,-35)(0)[t]{0}
\end{picture}
}}
\\ \hline
%
$D_{n+1}^{(2)}$ &
\raisebox{-0.7cm}{\scalebox{0.8}{
\begin{picture}(350,62)(-10,-12)
\BText(0,0){1}
\LongArrow(10,0)(30,0)
\BText(40,0){2}
\LongArrow(50,0)(70,0)
\Text(85,0)[]{$\cdots$}
\LongArrow(95,0)(115,0)
\BText(125,0){n}
\LongArrow(135,0)(155,0)
\BText(165,0){0}
\LongArrow(175,0)(195,0)
\BBoxc(205,0)(13,13)
\Text(205,0)[]{\footnotesize$\overline{\mbox{n}}$}
\LongArrow(215,0)(235,0)
\Text(250,0)[]{$\cdots$}
\LongArrow(260,0)(280,0)
\BBoxc(290,0)(13,13)
\Text(290,0)[]{\footnotesize$\overline{\mbox{2}}$}
\LongArrow(300,0)(320,0)
\BBoxc(330,0)(13,13)
\Text(330,0)[]{\footnotesize$\overline{\mbox{1}}$}
\Text(165,35)[]{\footnotesize$\emptyset$}
\LongArrowArc(165,-433)(468,71,89)
\LongArrowArc(165,-433)(468,91,109)
\PText(20,2)(0)[b]{1}
\PText(60,2)(0)[b]{2}
\PText(105,2)(0)[b]{n-1}
\PText(145,2)(0)[b]{n}
\PText(185,2)(0)[b]{n}
\PText(225,2)(0)[b]{n-1}
\PText(270,2)(0)[b]{2}
\PText(310,2)(0)[b]{1}
\PText(248,30)(0)[b]{0}
\PText(82,30)(0)[b]{0}
\end{picture}
}}
\\ \hline
\end{tabular}\vspace{4mm}
\caption{\label{table crystals}Level one perfect crystals}
\end{table}

The crystal $B(2\La_1)\oplus B(0)$ is a level one perfect 
crystal for $C_n^{(1)}$ ($n\ge 2$)~\cite{KKM}. 
The crystal graph corresponding to the 
integrable highest weight module $V(\La_1)$ of $U_q(C_n)$ is given
by~\cite[(4.2.4)]{KN}
\begin{center}
\raisebox{-0.2cm}{\scalebox{0.8}{
\begin{picture}(310,20)(-10,-10)
\BText(0,0){1}
\LongArrow(10,0)(30,0)
\BText(40,0){2}
\LongArrow(50,0)(70,0)
\Text(85,0)[]{$\cdots$}
\LongArrow(95,0)(115,0)
\BText(125,0){n}
\LongArrow(135,0)(155,0)
\BBoxc(165,0)(13,13)
\Text(165,0)[]{\footnotesize$\overline{\mbox{n}}$}
\LongArrow(175,0)(195,0)
\Text(210,0)[]{$\cdots$}
\LongArrow(220,0)(240,0)
\BBoxc(250,0)(13,13)
\Text(250,0)[]{\footnotesize$\overline{\mbox{2}}$}
\LongArrow(260,0)(280,0)
\BBoxc(290,0)(13,13)
\Text(290,0)[]{\footnotesize$\overline{\mbox{1}}$}
\PText(20,2)(0)[b]{1}
\PText(60,2)(0)[b]{2}
\PText(105,2)(0)[b]{n-1}
\PText(145,2)(0)[b]{n}
\PText(185,2)(0)[b]{n-1}
\PText(230,2)(0)[b]{2}
\PText(270,2)(0)[b]{1}
\end{picture}}}.
\end{center}
The crystal $B(2\La_1)$ is the connected component of
$B(\La_1)\otimes B(\La_1)$ containing $u_{\La_1}\otimes u_{\La_1}$
(see~\cite[Section 4.4]{KN}) which fixes the action of $e_i$ and $f_i$ 
for $1\le i\le n$. The edges in $B(2\La_1)\oplus B(0)$
corresponding to $f_0$ are given by~\cite{KKM}
\begin{center}
\begin{tabular}{ll}
\raisebox{-0.15cm}{\scalebox{0.8}{
\begin{picture}(94,22)(-2,-10)
\BBoxc(13,0)(26,13)
\Line(13,-6.5)(13,6.5)
\Text(6,0)[]{\footnotesize$\mbox{i}$}
\Text(19,0)[]{\footnotesize$\overline{\mbox{1}}$}
\LongArrow(30,0)(60,0)
\BBoxc(77,0)(26,13)
\Line(77,-6.5)(77,6.5)
\Text(70,0)[]{\footnotesize$\mbox{1}$}
\Text(83,0)[]{\footnotesize$\mbox{i}$}
\PText(45,2)(0)[b]{0}
\end{picture}}}
& for $i\neq 1,\overline{1}$ \\
\raisebox{-0.15cm}{\scalebox{0.8}{
\begin{picture}(94,22)(-2,-10)
\BBoxc(13,0)(26,13)
\Line(13,-6.5)(13,6.5)
\Text(6,0)[]{\footnotesize$\overline{\mbox{1}}$}
\Text(19,0)[]{\footnotesize$\overline{\mbox{1}}$}
\LongArrow(30,0)(60,0)
\Text(70,0)[]{\footnotesize$\emptyset$}
\PText(45,2)(0)[b]{0}
\end{picture}}}
&\\
\raisebox{-0.15cm}{\scalebox{0.8}{
\begin{picture}(94,22)(-2,-10)
\Text(19,0)[]{\footnotesize$\emptyset$}
\LongArrow(30,0)(60,0)
\BBoxc(77,0)(26,13)
\Line(77,-6.5)(77,6.5)
\Text(70,0)[]{\footnotesize$\mbox{1}$}
\Text(83,0)[]{\footnotesize$\mbox{1}$}
\PText(45,2)(0)[b]{0}
\end{picture}}}\,.
& \end{tabular}
\end{center}
There are the following strings of length greater than one
\begin{subequations}
\begin{align}
&\raisebox{-0.15cm}{\scalebox{0.8}{
\begin{picture}(240,22)(-30,-10)
\BBoxc(0,0)(26,13)
\Line(0,-6.5)(0,6.5)
\Text(-7,0)[]{\footnotesize$\mbox{k}$}
\Text(6,0)[]{\footnotesize$\mbox{k}$}
\LongArrow(30,0)(60,0)
\BBoxc(96.5,0)(39,13)
\Line(90,-6.5)(90,6.5)
\Text(83,0)[]{\footnotesize$\mbox{k}$}
\Text(103,0)[]{\footnotesize$\mbox{k+1}$}
\LongArrow(120,0)(150,0)
\BBoxc(180,0)(52,13)
\Line(180,-6.5)(180,6.5)
\Text(167,0)[]{\footnotesize$\mbox{k+1}$}
\Text(193,0)[]{\footnotesize$\mbox{k+1}$}
\PText(45,2)(0)[b]{k}
\PText(135,2)(0)[b]{k}
\end{picture}}}
&&  \text{for $1\leq k<n$} \notag \\
&\raisebox{-0.15cm}{\scalebox{0.8}{
\begin{picture}(240,22)(-30,-10)
\BBoxc(6.5,0)(39,13)
\Line(0,-6.5)(0,6.5)
\Text(-7,0)[]{\footnotesize$\mbox{k}$}
\Text(13,0)[]{\footnotesize$\overline{\mbox{k+1}}$}
\LongArrow(30,0)(60,0)
\BBoxc(90,0)(26,13)
\Line(90,-6.5)(90,6.5)
\Text(83,0)[]{\footnotesize$\mbox{k}$}
\Text(96,0)[]{\footnotesize$\overline{\mbox{k}}$}
\LongArrow(120,0)(150,0)
\BBoxc(173.5,0)(39,13)
\Line(180,-6.5)(180,6.5)
\Text(167,0)[]{\footnotesize$\mbox{k+1}$}
\Text(186,0)[]{\footnotesize$\overline{\mbox{k}}$}
\PText(45,2)(0)[b]{k}
\PText(135,2)(0)[b]{k}
\end{picture}}}
&& \text{for $1\leq k<n$} \label{string gen}\\
&\raisebox{-0.15cm}{\scalebox{0.8}{
\begin{picture}(240,22)(-30,-10)
\BBoxc(0,0)(52,13)
\Line(0,-6.5)(0,6.5)
\Text(-13,0)[]{\footnotesize$\overline{\mbox{k+1}}$}
\Text(13,0)[]{\footnotesize$\overline{\mbox{k+1}}$}
\LongArrow(30,0)(60,0)
\BBoxc(83.5,0)(39,13)
\Line(90,-6.5)(90,6.5)
\Text(77,0)[]{\footnotesize$\overline{\mbox{k+1}}$}
\Text(96,0)[]{\footnotesize$\overline{\mbox{k}}$}
\LongArrow(120,0)(150,0)
\BBoxc(180,0)(26,13)
\Line(180,-6.5)(180,6.5)
\Text(173,0)[]{\footnotesize$\overline{\mbox{k}}$}
\Text(186,0)[]{\footnotesize$\overline{\mbox{k}}$}
\PText(45,2)(0)[b]{k}
\PText(135,2)(0)[b]{k}
\end{picture}}}
&& \text{for $1\leq k<n$} \notag \\[2mm]
&\raisebox{-0.15cm}{\scalebox{0.8}{
\begin{picture}(240,22)(-30,-10)
\BBoxc(0,0)(26,13)
\Line(0,-6.5)(0,6.5)
\Text(-7,0)[]{\footnotesize$\mbox{n}$}
\Text(6,0)[]{\footnotesize$\mbox{n}$}
\LongArrow(30,0)(60,0)
\BBoxc(90,0)(26,13)
\Line(90,-6.5)(90,6.5)
\Text(83,0)[]{\footnotesize$\mbox{n}$}
\Text(96,0)[]{\footnotesize$\overline{\mbox{n}}$}
\LongArrow(120,0)(150,0)
\BBoxc(180,0)(26,13)
\Line(180,-6.5)(180,6.5)
\Text(173,0)[]{\footnotesize$\overline{\mbox{n}}$}
\Text(186,0)[]{\footnotesize$\overline{\mbox{n}}$}
\PText(45,2)(0)[b]{n}
\PText(135,2)(0)[b]{n}
\end{picture}}}
&& \label{string n}\\[2mm]
&\raisebox{-0.15cm}{\scalebox{0.8}{
\begin{picture}(240,22)(-30,-10)
\BBoxc(0,0)(26,13)
\Line(0,-6.5)(0,6.5)
\Text(-7,0)[]{\footnotesize$\overline{\mbox{1}}$}
\Text(6,0)[]{\footnotesize$\overline{\mbox{1}}$}
\LongArrow(30,0)(60,0)
\Text(90,0)[]{\footnotesize$\emptyset$}
\LongArrow(120,0)(150,0)
\BBoxc(180,0)(26,13)
\Line(180,-6.5)(180,6.5)
\Text(173,0)[]{\footnotesize$\mbox{1}$}
\Text(186,0)[]{\footnotesize$\mbox{1}$}
\PText(45,2)(0)[b]{0}
\PText(135,2)(0)[b]{0}
\end{picture}}}.\notag
&&
\end{align}
\end{subequations}
Note that none of the elements satisfies \eqref{cond}.

For type $A_{2n-1}^{(2)}$ the crystal $B(2\La_1)$ is perfect of 
level 2 \cite[Sec. 1.6 and 6.7]{KMN2}.
The elements are given by
$\raisebox{-0.15cm}{\scalebox{0.8}{\begin{picture}(29,20)(-2,-10)
\BBoxc(13,0)(26,13)
\Line(13,-6.5)(13,6.5)
\Text(6,0)[]{\footnotesize$\mbox{x}$}
\Text(19,0)[]{\footnotesize$\mbox{y}$}
\end{picture}}}$ 
with $x\le y$ and 
$x,y \in \{1<2<\cdots<n<\bar{n}<\cdots <\bar{2}<\bar{1} \}$.
The action of $f_i$ for $i=1,2,\ldots,n$ is the same as for
the above $C_n^{(1)}$ crystal of level one, 
and $f_0=\sigma\circ f_1\circ \sigma$
where $\sigma$ is the involution that exchanges 1 and $\bar{1}$ (with
appropriate reorderings).

The strings of length greater than one are the same as in
\eqref{string gen} and \eqref{string n}. In addition there are
the following $0$-strings of length 2
\begin{align}
&\raisebox{-0.15cm}{\scalebox{0.8}{
\begin{picture}(240,22)(-30,-10)
\BBoxc(0,0)(26,13)
\Line(0,-6.5)(0,6.5)
\Text(-7,0)[]{\footnotesize$\overline{\mbox{1}}$}
\Text(6,0)[]{\footnotesize$\overline{\mbox{1}}$}
\LongArrow(30,0)(60,0)
\BBoxc(90,0)(26,13)
\Line(90,-6.5)(90,6.5)
\Text(83,0)[]{\footnotesize$\mbox{2}$}
\Text(96,0)[]{\footnotesize$\overline{\mbox{1}}$}
\LongArrow(120,0)(150,0)
\BBoxc(180,0)(26,13)
\Line(180,-6.5)(180,6.5)
\Text(173,0)[]{\footnotesize$\mbox{2}$}
\Text(186,0)[]{\footnotesize$\mbox{2}$}
\PText(45,2)(0)[b]{0}
\PText(135,2)(0)[b]{0}
\end{picture}}} 
&& \notag \\
&\raisebox{-0.15cm}{\scalebox{0.8}{
\begin{picture}(240,22)(-30,-10)
\BBoxc(0,0)(26,13)
\Line(0,-6.5)(0,6.5)
\Text(-7,0)[]{\footnotesize$\overline{\mbox{2}}$}
\Text(6,0)[]{\footnotesize$\overline{\mbox{1}}$}
\LongArrow(30,0)(60,0)
\BBoxc(90,0)(26,13)
\Line(90,-6.5)(90,6.5)
\Text(83,0)[]{\footnotesize$\mbox{1}$}
\Text(96,0)[]{\footnotesize$\overline{\mbox{1}}$}
\LongArrow(120,0)(150,0)
\BBoxc(180,0)(26,13)
\Line(180,-6.5)(180,6.5)
\Text(173,0)[]{\footnotesize$\mbox{1}$}
\Text(186,0)[]{\footnotesize$\mbox{2}$}
\PText(45,2)(0)[b]{0}
\PText(135,2)(0)[b]{0}
\end{picture}}} 
&& \label{string 0}\\
&\raisebox{-0.15cm}{\scalebox{0.8}{
\begin{picture}(240,22)(-30,-10)
\BBoxc(0,0)(26,13)
\Line(0,-6.5)(0,6.5)
\Text(-7,0)[]{\footnotesize$\overline{\mbox{2}}$}
\Text(6,0)[]{\footnotesize$\overline{\mbox{2}}$}
\LongArrow(30,0)(60,0)
\BBoxc(90,0)(26,13)
\Line(90,-6.5)(90,6.5)
\Text(83,0)[]{\footnotesize$\mbox{1}$}
\Text(96,0)[]{\footnotesize$\overline{\mbox{2}}$}
\LongArrow(120,0)(150,0)
\BBoxc(180,0)(26,13)
\Line(180,-6.5)(180,6.5)
\Text(173,0)[]{\footnotesize$\mbox{1}$}
\Text(186,0)[]{\footnotesize$\mbox{1}$}
\PText(45,2)(0)[b]{0}
\PText(135,2)(0)[b]{0}
\end{picture}}}. 
&& \notag
\end{align}
The only elements fulfilling \eqref{cond} are
$\raisebox{-0.15cm}{\scalebox{0.8}{\begin{picture}(29,20)(-2,-10)
\BBoxc(13,0)(26,13)
\Line(13,-6.5)(13,6.5)
\Text(6,0)[]{\footnotesize$\mbox{1}$}
\Text(19,0)[]{\footnotesize$\overline{\mbox{1}}$}
\end{picture}}}$,
$\raisebox{-0.15cm}{\scalebox{0.8}{\begin{picture}(29,20)(-2,-10)
\BBoxc(13,0)(26,13)
\Line(13,-6.5)(13,6.5)
\Text(6,0)[]{\footnotesize$\mbox{1}$}
\Text(19,0)[]{\footnotesize$\mbox{2}$}
\end{picture}}}$,
$\raisebox{-0.15cm}{\scalebox{0.8}{\begin{picture}(29,20)(-2,-10)
\BBoxc(13,0)(26,13)
\Line(13,-6.5)(13,6.5)
\Text(6,0)[]{\footnotesize$\mbox{2}$}
\Text(19,0)[]{\footnotesize$\overline{\mbox{1}}$}
\end{picture}}}$, and
$\raisebox{-0.15cm}{\scalebox{0.8}{\begin{picture}(29,20)(-2,-10)
\BBoxc(13,0)(26,13)
\Line(13,-6.5)(13,6.5)
\Text(6,0)[]{\footnotesize$\mbox{2}$}
\Text(19,0)[]{\footnotesize$\mbox{2}$}
\end{picture}}}$
which belong to a $0$-string and a $1$-string of length two.
It can be checked that setting $v(b)=0$ for $b$ one of these four 
elements guarantees the involutive condition of $v$.

For type $B_n^{(1)}$ the crystal $B(2\La_1)$ is perfect of 
level 2 \cite[Sec. 1.7 and 6.8]{KMN2}. It consists of the elements
$\raisebox{-0.15cm}{\scalebox{0.8}{\begin{picture}(29,20)(-2,-10)
\BBoxc(13,0)(26,13)
\Line(13,-6.5)(13,6.5)
\Text(6,0)[]{\footnotesize$\mbox{x}$}
\Text(19,0)[]{\footnotesize$\mbox{y}$}
\end{picture}}}$ 
with $x\le y$ and 
$x,y\in \{1<\cdots<n<0<\ba{n}<\cdots<\ba{1}\}$;
$x=y=0$ is excluded.
The action of $f_i$ for $i=1,2,\ldots,n$ is given by
the tensor product rule using the action on the level 1 crystal
of $B_n^{(1)}$ as given in Table \ref{table crystals}, and 
$f_0=\sigma\circ f_1\circ \sigma$ where $\sigma$ is the involution 
that exchanges 1 and $\bar{1}$ (with appropriate reorderings).

The strings of length greater than one are those of
equations \eqref{string gen} and \eqref{string 0} and in addition
the following $n$-string of length four
\begin{equation}\label{string n length 4}
\raisebox{-0.15cm}{\scalebox{0.8}{
\begin{picture}(330,22)(-30,-10)
\BBoxc(0,0)(26,13)
\Line(0,-6.5)(0,6.5)
\Text(-7,0)[]{\footnotesize$\mbox{n}$}
\Text(6,0)[]{\footnotesize$\mbox{n}$}
\LongArrow(20,0)(50,0)
\BBoxc(70,0)(26,13)
\Line(70,-6.5)(70,6.5)
\Text(63,0)[]{\footnotesize$\mbox{n}$}
\Text(76,0)[]{\footnotesize$\mbox{0}$}
\LongArrow(90,0)(120,0)
\BBoxc(140,0)(26,13)
\Line(140,-6.5)(140,6.5)
\Text(133,0)[]{\footnotesize$\mbox{n}$}
\Text(146,0)[]{\footnotesize$\overline{\mbox{n}}$}
\LongArrow(160,0)(190,0)
\BBoxc(210,0)(26,13)
\Line(210,-6.5)(210,6.5)
\Text(203,0)[]{\footnotesize$\mbox{0}$}
\Text(216,0)[]{\footnotesize$\overline{\mbox{n}}$}
\LongArrow(230,0)(260,0)
\BBoxc(280,0)(26,13)
\Line(280,-6.5)(280,6.5)
\Text(273,0)[]{\footnotesize$\overline{\mbox{n}}$}
\Text(286,0)[]{\footnotesize$\overline{\mbox{n}}$}
\PText(35,2)(0)[b]{n}
\PText(105,2)(0)[b]{n}
\PText(175,2)(0)[b]{n}
\PText(245,2)(0)[b]{n}
\end{picture}}}.
\end{equation}
The same four elements as for $A_{2n-1}^{(2)}$ satisfy \eqref{cond}
and again setting $v(b)=0$ for these ensures the involutive property of $v$.

For type $D_n^{(1)}$ the crystal $B(2\La_1)$ is perfect of 
level 2 \cite[Sec. 1.8 and 6.9]{KMN2}. It consists of the elements
$\raisebox{-0.15cm}{\scalebox{0.8}{\begin{picture}(29,20)(-2,-10)
\BBoxc(13,0)(26,13)
\Line(13,-6.5)(13,6.5)
\Text(6,0)[]{\footnotesize$\mbox{x}$}
\Text(19,0)[]{\footnotesize$\mbox{y}$}
\end{picture}}}$ 
with $x\le y$ and 
$x,y\in \{1<2<\cdots<n,\ba{n}<\cdots<\ba{1}\}$,
the cases $x=n, y=\ba{n}$ and $x=\ba{n}, y=n$ being excluded.
The action of $f_i$ for $i=1,2,\ldots,n$ is given by
the tensor product rule using the action on the level 1 crystal
of $D_n^{(1)}$ as given in Table \ref{table crystals}, and 
$f_0=\sigma\circ f_1\circ \sigma$ where $\sigma$ is the involution 
that exchanges 1 and $\bar{1}$ (with appropriate reorderings).

Again the strings of length greater than one are the same as in
equations \eqref{string gen} and \eqref{string 0} plus the
following $n$-strings
\begin{align*}
&\raisebox{-0.15cm}{\scalebox{0.8}{
\begin{picture}(240,22)(-30,-10)
\BBoxc(0,0)(26,13)
\Line(0,-6.5)(0,6.5)
\Text(-7,0)[]{\footnotesize$\mbox{n}$}
\Text(6,0)[]{\footnotesize$\mbox{n}$}
\LongArrow(30,0)(60,0)
\BBoxc(96.5,0)(39,13)
\Line(90,-6.5)(90,6.5)
\Text(83,0)[]{\footnotesize$\mbox{n}$}
\Text(103,0)[]{\footnotesize$\overline{\mbox{n-1}}$}
\LongArrow(120,0)(150,0)
\BBoxc(180,0)(52,13)
\Line(180,-6.5)(180,6.5)
\Text(167,0)[]{\footnotesize$\overline{\mbox{n-1}}$}
\Text(193,0)[]{\footnotesize$\overline{\mbox{n-1}}$}
\PText(45,2)(0)[b]{n}
\PText(135,2)(0)[b]{n}
\end{picture}}}\\
&\raisebox{-0.15cm}{\scalebox{0.8}{
\begin{picture}(240,22)(-30,-10)
\BBoxc(0,0)(52,13)
\Line(0,-6.5)(0,6.5)
\Text(-13,0)[]{\footnotesize$\mbox{n-1}$}
\Text(13,0)[]{\footnotesize$\mbox{n-1}$}
\LongArrow(30,0)(60,0)
\BBoxc(83.5,0)(39,13)
\Line(90,-6.5)(90,6.5)
\Text(77,0)[]{\footnotesize$\mbox{n-1}$}
\Text(96,0)[]{\footnotesize$\overline{\mbox{n}}$}
\LongArrow(120,0)(150,0)
\BBoxc(180,0)(26,13)
\Line(180,-6.5)(180,6.5)
\Text(173,0)[]{\footnotesize$\overline{\mbox{n}}$}
\Text(186,0)[]{\footnotesize$\overline{\mbox{n}}$}
\PText(45,2)(0)[b]{n}
\PText(135,2)(0)[b]{n}
\end{picture}}}\\
&\raisebox{-0.15cm}{\scalebox{0.8}{
\begin{picture}(240,22)(-30,-10)
\BBoxc(-6.5,0)(39,13)
\Line(0,-6.5)(0,6.5)
\Text(-13,0)[]{\footnotesize$\mbox{n-1}$}
\Text(6.5,0)[]{\footnotesize$\mbox{n}$}
\LongArrow(30,0)(60,0)
\BBoxc(90,0)(52,13)
\Line(90,-6.5)(90,6.5)
\Text(77,0)[]{\footnotesize$\mbox{n-1}$}
\Text(103,0)[]{\footnotesize$\overline{\mbox{n-1}}$}
\LongArrow(120,0)(150,0)
\BBoxc(186.5,0)(39,13)
\Line(180,-6.5)(180,6.5)
\Text(173,0)[]{\footnotesize$\overline{\mbox{n}}$}
\Text(193,0)[]{\footnotesize$\overline{\mbox{n-1}}$}
\PText(45,2)(0)[b]{n}
\PText(135,2)(0)[b]{n}
\end{picture}}}\; .
\end{align*}

In addition to the four elements
$\raisebox{-0.15cm}{\scalebox{0.8}{\begin{picture}(29,20)(-2,-10)
\BBoxc(13,0)(26,13)
\Line(13,-6.5)(13,6.5)
\Text(6,0)[]{\footnotesize$\mbox{1}$}
\Text(19,0)[]{\footnotesize$\overline{\mbox{1}}$}
\end{picture}}}$,
$\raisebox{-0.15cm}{\scalebox{0.8}{\begin{picture}(29,20)(-2,-10)
\BBoxc(13,0)(26,13)
\Line(13,-6.5)(13,6.5)
\Text(6,0)[]{\footnotesize$\mbox{1}$}
\Text(19,0)[]{\footnotesize$\mbox{2}$}
\end{picture}}}$,
$\raisebox{-0.15cm}{\scalebox{0.8}{\begin{picture}(29,20)(-2,-10)
\BBoxc(13,0)(26,13)
\Line(13,-6.5)(13,6.5)
\Text(6,0)[]{\footnotesize$\mbox{2}$}
\Text(19,0)[]{\footnotesize$\overline{\mbox{1}}$}
\end{picture}}}$, and
$\raisebox{-0.15cm}{\scalebox{0.8}{\begin{picture}(29,20)(-2,-10)
\BBoxc(13,0)(26,13)
\Line(13,-6.5)(13,6.5)
\Text(6,0)[]{\footnotesize$\mbox{2}$}
\Text(19,0)[]{\footnotesize$\mbox{2}$}
\end{picture}}}$
also the elements
$\raisebox{-0.15cm}{\scalebox{0.8}{\begin{picture}(60,20)(-2,-10)
\BBoxc(26,0)(52,13)
\Line(26,-6.5)(26,6.5)
\Text(12,0)[]{\footnotesize$\overline{\mbox{n-1}}$}
\Text(38,0)[]{\footnotesize$\overline{\mbox{n-1}}$}
\end{picture}}}$,
$\raisebox{-0.15cm}{\scalebox{0.8}{\begin{picture}(60,20)(-2,-10)
\BBoxc(26,0)(52,13)
\Line(26,-6.5)(26,6.5)
\Text(12,0)[]{\footnotesize$\mbox{n-1}$}
\Text(38,0)[]{\footnotesize$\overline{\mbox{n-1}}$}
\end{picture}}}$,
$\raisebox{-0.15cm}{\scalebox{0.8}{\begin{picture}(44,20)(-2,-10)
\BBoxc(19.5,0)(39,13)
\Line(13,-6.5)(13,6.5)
\Text(6,0)[]{\footnotesize$\mbox{n}$}
\Text(25,0)[]{\footnotesize$\overline{\mbox{n-1}}$}
\end{picture}}}$, and
$\raisebox{-0.15cm}{\scalebox{0.8}{\begin{picture}(44,20)(-2,-10)
\BBoxc(19.5,0)(39,13)
\Line(13,-6.5)(13,6.5)
\Text(6,0)[]{\footnotesize$\overline{\mbox{n}}$}
\Text(25,0)[]{\footnotesize$\overline{\mbox{n-1}}$}
\end{picture}}}$
satisfy \eqref{cond}. The latter ones are contained in an 
$(n-1)$-string and an $n$-string. Setting $v(b)=0$ for the 
first four elements and $v(b)=n$ for the last four elements ensures 
the involutive property of $v$.

The crystal $B(0)\oplus B(\Lambda_1)\oplus B(2\Lambda_1)$
is a level 2 perfect crystal for $D_{n+1}^{(2)}$
\cite[Sections 1.9 and 6.10]{KMN2}.
The elements of this crystal are $\emptyset$, $\boxf{x}$, and
$\raisebox{-0.15cm}{\scalebox{0.8}{\begin{picture}(29,20)(-2,-10)
\BBoxc(13,0)(26,13)
\Line(13,-6.5)(13,6.5)
\Text(6,0)[]{\footnotesize$\mbox{x}$}
\Text(19,0)[]{\footnotesize$\mbox{y}$}
\end{picture}}}$ 
with $x,y\in \{1<2<\cdots<n<0<\ba{n}<\cdots<\ba{1}\}$ and
$x\le y$; $x=y=0$ is excluded.
The action of $f_i$ for $i=1,2,\ldots,n$ is given by
the tensor product rule using the action on the level 1 crystal
of $D_{n+1}^{(2)}$ as given in Table \ref{table crystals}, and the 
action of $f_0$ is given by
\begin{equation}\label{f_0 action}
\begin{split}
\raisebox{-0.15cm}{\scalebox{0.8}{
\begin{picture}(94,22)(-2,-10)
\Text(19,0)[]{\footnotesize$\emptyset$}
\LongArrow(30,0)(60,0)
\BBoxc(70,0)(13,13)
\Text(70,0)[]{\footnotesize$\mbox{1}$}
\PText(45,2)(0)[b]{0}
\end{picture}}}
&\\
\raisebox{-0.15cm}{\scalebox{0.8}{
\begin{picture}(94,22)(-2,-10)
\BBoxc(19.5,0)(13,13)
\Text(19,0)[]{\footnotesize$\mbox{x}$}
\LongArrow(30,0)(60,0)
\BBoxc(77,0)(26,13)
\Line(77,-6.5)(77,6.5)
\Text(70,0)[]{\footnotesize$\mbox{1}$}
\Text(83,0)[]{\footnotesize$\mbox{x}$}
\PText(45,2)(0)[b]{0}
\end{picture}}}
&\qquad \text{for $x\neq \overline{1}$}\\
\raisebox{-0.15cm}{\scalebox{0.8}{
\begin{picture}(94,22)(-2,-10)
\BBoxc(19.5,0)(13,13)
\Text(19,0)[]{\footnotesize$\overline{\mbox{1}}$}
\LongArrow(30,0)(60,0)
\Text(70,0)[]{\footnotesize$\emptyset$}
\PText(45,2)(0)[b]{0}
\end{picture}}}
&\\
\raisebox{-0.15cm}{\scalebox{0.8}{
\begin{picture}(94,22)(-2,-10)
\BBoxc(13,0)(26,13)
\Line(13,-6.5)(13,6.5)
\Text(6,0)[]{\footnotesize$\mbox{x}$}
\Text(19,0)[]{\footnotesize$\overline{\mbox{1}}$}
\LongArrow(30,0)(60,0)
\BBoxc(70,0)(13,13)
\Text(70,0)[]{\footnotesize$\mbox{x}$}
\PText(45,2)(0)[b]{0}
\end{picture}}}
& \qquad \text{for $x\neq 1$}
\end{split}
\end{equation}
and undefined otherwise.

The strings of length greater than one are given by
\eqref{string gen}, \eqref{string n length 4} and
\begin{align}
\label{string 0 length 4}
&\raisebox{-0.15cm}{\scalebox{0.8}{
\begin{picture}(280,22)(-30,-10)
\BBoxc(0,0)(26,13)
\Line(0,-6.5)(0,6.5)
\Text(-7,0)[]{\footnotesize$\overline{\mbox{1}}$}
\Text(6,0)[]{\footnotesize$\overline{\mbox{1}}$}
\LongArrow(20,0)(50,0)
\BBoxc(65,0)(13,13)
\Text(66,)[]{\footnotesize$\overline{\mbox{1}}$}
\LongArrow(80,0)(110,0)
\Text(125,0)[]{\footnotesize$\emptyset$}
\LongArrow(140,0)(170,0)
\BBoxc(185,0)(13,13)
\Text(186,0)[]{\footnotesize$\mbox{1}$}
\LongArrow(200,0)(230,0)
\BBoxc(250,0)(26,13)
\Line(250,-6.5)(250,6.5)
\Text(243,0)[]{\footnotesize$\mbox{1}$}
\Text(256,0)[]{\footnotesize$\mbox{1}$}
\PText(35,2)(0)[b]{0}
\PText(95,2)(0)[b]{0}
\PText(155,2)(0)[b]{0}
\PText(215,2)(0)[b]{0}
\end{picture}}}
\\[2mm]
&\raisebox{-0.15cm}{\scalebox{0.8}{
\begin{picture}(280,22)(-30,-10)
\BBoxc(65,0)(13,13)
\Text(66,0)[]{\footnotesize$\mbox{n}$}
\LongArrow(80,0)(110,0)
\BBoxc(125,0)(13,13)
\Text(126,0)[]{\footnotesize$\mbox{0}$}
\LongArrow(140,0)(170,0)
\BBoxc(185,0)(13,13)
\Text(186,0)[]{\footnotesize$\overline{\mbox{n}}$}
\PText(95,2)(0)[b]{n}
\PText(155,2)(0)[b]{n}
\end{picture}}}
\end{align}
There are no elements with property \eqref{cond}.

The crystal $B(0)\oplus B(\Lambda_1)\oplus B(2\Lambda_1)$
is a level 2 perfect crystal for $A_{2n}^{(2)}$
\cite[Sec. 1.10 and 6.11]{KMN2}.
The elements of this crystal are $\emptyset$, $\boxf{x}$, and
$\raisebox{-0.15cm}{\scalebox{0.8}{\begin{picture}(29,20)(-2,-10)
\BBoxc(13,0)(26,13)
\Line(13,-6.5)(13,6.5)
\Text(6,0)[]{\footnotesize$\mbox{x}$}
\Text(19,0)[]{\footnotesize$\mbox{y}$}
\end{picture}}}$ 
with $x,y\in \{1<2<\cdots<n<\ba{n}<\cdots<\ba{1}\}$ and $x\le y$.
The action of $f_i$ for $i=1,2,\ldots,n$ is given by
the tensor product rule using the action on the level 1 crystal 
of $A_{2n}^{(2)}$ as given in Table \ref{table crystals}, and the action of
$f_0$ is the same as in \eqref{f_0 action}.

The strings of length greater than one are as in \eqref{string gen} for
$n\ge 2$, \eqref{string n} and \eqref{string 0 length 4}.
Again there are no elements with property \eqref{cond}.
\end{proof}

\begin{rem} \label{northeast} Suppose $\gggg$ is of type
$A^{(1)}_{n-1}$ in Lemma \ref{lev one}.
The function $v$ amounts to a canonical choice of 
a simple root $i$ among those such that the given
element admits $e_i$.  Consider $b\in \BBB(\La_r)$
such that $b \not= u_{\La_r}$.  In addition to the
realization of the crystal $\BBB(\La_r)$ by the
space of homogeneous paths using the crystal given
in the proof of Lemma \ref{lev one}, one may also
consider the realization in \cite{DJKMO} by
$n$-regular partitions.  Suppose $\la$ is the
partition corresponding to $b$.  Then
up to the Dynkin diagram automorphism that sends
$r+i$ to $r-i$ modulo $n$, the choice of violation
$v$ corresponds to the corner cell of $\la$
that is in the rightmost column of $\la$.
This choice of corner cell is used in \cite{Lak}
to define the smallest Demazure crystal of
$\BBB(\La_r)$ containing $b$.
\end{rem}

\section{Inhomogeneous paths}
\label{sec inhomo paths}

\begin{thm} \label{pi} 
Let $\gggg$ be as in Theorem \ref{gen boson}, and $B$, $\La$, and $B_0$ 
be as in \eqref{wt tensor iso}.  
Suppose in addition that
for all $1\le j\le L$ and $b\in B_j$, if
$b\otimes b_0\mapsto b_0'\otimes b'$ under the local isomorphism
$B_j \otimes B_0 \rightarrow B_0 \otimes B_j$ and
$e_0(b \otimes b_0) = e_0(b) \otimes b_0$ then
$e_0(b_0'\otimes b') = e_0(b_0') \otimes b'$.
Then 
\begin{equation} \label{pi formula}
  \Pi\,(\ch(B \otimes u_\La)) = \ch(B \otimes \BBB(\La)).
\end{equation}
\end{thm}
\begin{proof} It is enough to verify the hypotheses of
Theorem \ref{gen boson}, applied to
$B \otimes u_\La \subset B \otimes \BBB(\La)$.
$B \otimes \BBB(\La)$ is isomorphic to a direct sum of irreducible integrable
highest weight modules by Theorem \ref{level tensor}.
$B\otimes u_\La$ is obviously closed under the $e_i$.
It follows from \cite[Lemma 1]{KMOU} that
$B \otimes u_\La$ generates $B \otimes \BBB(\La)$.
To check the third condition of Theorem \ref{gen boson},
let $b\in B$ and $i\in I$ be such that
$\epsilon_i(b \otimes u_\La) > 0$.  Then
$\epsilon_i(b) > \phi_i(u_\La) = \inner{\ack_i}{\La}$.
This implies that the $i$-string of $b \otimes u_\La$ inside
$B \otimes \BBB(\La)$, consists of vectors of the form $b' \otimes u_\La$
where $b' \in B$.

Finally, Lemma \ref{energy e} with $B$ replaced by $B\otimes B_0$
guarantees that the
affine weight function on $B \otimes \BBB(\La)$
determined by its value on highest weight
vectors, agrees on the subset $B\otimes u_\La$ with the function
$\wt(b) = \af(\wt'(b)) - E_{B,B_0}(b \otimes b_0) a_0^{-1} \delta$
where $\wt':B\rightarrow \Pcl$ is the original weight function.
\end{proof}

\begin{rem} 
Observe that even without the extra hypothesis on the action
of $e_0$ in Theorem \ref{pi}, one obtains a bosonic formula.  
The extra condition is only needed to show that
the energy function $b\mapsto E_{B,B_0}(b\otimes b_0)$
gives rise to the correct affine weight for all elements of the
form $b\otimes u_\La$ and not just on the highest weight vectors.
Perhaps this extra condition is always a consequence of the other hypotheses.
\end{rem}

Now the formula \eqref{pi formula} is written more explicitly.
Let $m\in \Z$ and $\La,\La'\in \af(\Pcl^+)$ be of level $\ell$.
A formula equivalent to \eqref{pi formula} is obtained by
taking the coefficient of $\ch \VVV(\La'-m a_0^{-1}\delta)$ on both sides:
\begin{equation*}
  [q^m] \KP(B,\La,\La',B_0)(q) = \sum_{(w,b)\in \SSS} (-1)^w
\end{equation*}
where $\SSS$ is the set of pairs $(w,b)\in W \times B$ such that
\begin{equation} \label{wt cond 2}
  w^{-1}(\La'+\rho)-m a_0^{-1}\delta -\rho = \wt(b \otimes u_\La).
\end{equation}
Let $M$ be the sublattice of $\Pfin$ given by the
image under $\nu$ of the $\Z$-span of the orbit $\Wfin \thck$.
Let $T \subset GL(\hhh^*)$ be the group of translations
by the elements of $M$, where $t_\alpha \in T$ is translation
by $\alpha\in M$.  Then $W \cong T \rtimes \Wfin$ and
$r_0 = t_{\nu(\thck)} r_\theta$.  For $\alpha\in M$ and
$\La\in P$ of level $\ell$, one has \cite[(6.5.2)]{Kac}
\begin{equation} \label{trans act}
  t_\alpha(\La) = \La + \ell \alpha -
  ((\La|\alpha)+\tfrac{1}{2}\,|\alpha|^2 \ell)\delta.
\end{equation}
The action of $\tau\in \Wfin$ on the level $\ell$ weight $\La$
is given by
\begin{equation*}
  \tau(\La) = \tau(\overline{\La}+\ell\La_0)=
  	\tau(\overline{\La}) + \ell \La_0.
\end{equation*}
Now $\rho=\hck \La_0 + \rhofin$ where $\hck$ is the 
dual Coxeter number
and $\rhofin$ is the half-sum of the positive roots in $\gfin$.

Recall that $\Wfin$ leaves $\delta$ invariant.
In \eqref{wt cond 2} write $w = t_\alpha \tau $ where
$\tau\in \Wfin$ and $\alpha\in M$, obtaining
\begin{equation*}
\begin{split}
  \wt(b \otimes u_\La) =&
  \tau^{-1} t_{-\alpha} (\La'+\rho)-m a_0^{-1} \delta-\rho \\
  =& -m a_0^{-1} \delta-\rho+\tau^{-1}\bigl\{
  	\La'+\rho - (\ell+\hck)\alpha \\
  	& - \{(\La'+\rho|-\alpha)+\tfrac{1}{2}\,|\alpha|^2 (\ell+\hck)\}\delta
	\bigr\} \\
  =& \ell\La_0-\rhofin+\tau^{-1}(\overline{\La'}+\rhofin-(\ell+\hck)\alpha)
  \\ &+
 \bigl\{-m a_0^{-1}+(\overline{\La'}+\rhofin|\alpha)
	-\tfrac{1}{2}\,|\alpha|^2(\ell+\hck) \bigr\}\delta
\end{split}  	
\end{equation*}
Since both sides are weights of level $\ell$, by
equating coefficients of $\delta$ and projections into $\Pfin$, one obtains
the equivalent conditions
\begin{equation} \label{fin wt}
  \overline{\wt(b)} = - \overline{\La} -\rhofin +
  	\tau^{-1}(\overline{\La'}-(\ell+\hck)\alpha+\rhofin)
\end{equation}
and
\begin{equation} \label{delta wt}
  a_0^{-1} E(b) = a_0^{-1} m -(\overline{\La'}+\rhofin|\alpha)+
  \tfrac{1}{2}\,|\alpha|^2(\ell+\hck).
\end{equation}
Therefore one has the equality
\begin{equation} \label{level K poly}
  \KP(B,\La,\La',B_0)(q) =
 \sum_{\tau\in \Wfin} \sum_{\alpha\in M} \sum_{b\in B}
 	(-1)^\tau q^{E(b)+a_0(\overline{\La'}+\rhofin\mid\alpha)-
 	\tfrac{1}{2}a_0 |\alpha|^2(\ell+\hck)}
\end{equation}
where $b\in B$ satisfies
\begin{equation*}
  \wt(b)=-\overline{\La}-\rhofin+
	\tau^{-1}(\overline{\La'}-(\ell+\hck)\alpha+\rhofin).
\end{equation*}

\section{Type A}

\subsection{Conjecture of \cite{FLOT}}

For simplicity let us assume that $\gggg$ is of untwisted affine type,
where $a_0=1$ and $(\rhofin|\theta)=\hck-1$ \cite[Ex. 6.2]{Kac}.

Let $\La \in P$ be a weight of level $\ell$ but not necessarily
dominant.  Consider the weight $\La+\rho$.  If it is regular
(not fixed by any $w\in W$) then there is a unique $w\in W$
such that $w (\La+\rho)\in P^+$.  It follows from the definition of
$\Pi$ that
\begin{equation} \label{Bott}
  \Pi\,e^{\La} = \begin{cases}
    (-1)^w
	\ch \VVV(w(\La+\rho)-\rho) &
	\text{if $\La+\rho$ is $W$-regular and }\\
	& \text{$w(\La+\rho)\in P^+$} \\
    0 & \text{if $\La+\rho$ is not $W$-regular.}
  \end{cases}
\end{equation}
Then for all $i\in I$,
\begin{equation} \label{refl rho}
  -\Pi\,e^\La = \Pi\,e^{r_i(\La+\rho)-\rho}.
\end{equation}
Suppose $i\not=0$.  Then
\begin{equation*}
\begin{split}
  r_i(\La+\rho)-\rho &=
  (\ell+\hck)\La_0 + r_i(\overline{\La}+\rhofin)-(\hck\La_0+\rhofin) \\
  &= \ell \La_0 - \alpha_i + r_i(\overline{\La}).
\end{split}
\end{equation*}
For $i=0$, recall that
\begin{equation*}
	r_0 = t_{\nu(\thck)}r_\theta=t_{\theta} r_\theta =
	r_\theta t_{-\theta}.
\end{equation*}
Then
\begin{equation*}
\begin{split}
   t_{-\theta}(\La+\rho) &=
  \La+\rho-(\ell+\hck)\theta+
   \{(\La+\rho|\theta)-\tfrac{1}{2}|\theta|^2(\ell+\hck)\}\delta \\
   &= (\ell+\hck)\La_0 + \rhofin + \Lab -(\ell+\hck)\theta
   +\{(\Lab|\theta)-(1+\ell)\}\delta
\end{split}
\end{equation*}
and
\begin{equation*}
\begin{split}
  r_0 (\La+\rho)-\rho &=
  r_\theta\bigl\{ (\ell+\hck)\La_0 + \rhofin + \Lab -(\ell+\hck)\theta \\
  &+\{(\Lab|\theta)-(1+\ell)\}\delta\bigr\}-\rho \\
  &= (\ell+\hck)\La_0 + \rhofin
	-\inner{\thck}{\rhofin}\theta + r_\theta(\Lab)\\
&+(\ell+\hck)\theta+\{(\Lab|\theta)-(1+\ell)\}\delta-(\hck\La_0+\rhofin)\\
  &= \ell \La_0 +r_\theta(\Lab)+(\ell+1)\theta+
    \{(\Lab|\theta)-(1+\ell)\}\delta.
\end{split}
\end{equation*}

Now let $\gggg$ be of type $A^{(1)}_{n-1}$.
Let $\Pfin$ be identified with the subspace of $\Z^n$
given by vectors with sum zero.

For $\alpha\in \Pfin$ define the Demazure operator
$\Pifin$ to be the linear operator on $\Z[\Pfin]$ such that
\begin{equation*}
s_\alpha := \Pifin(e^\alpha) =
\Jfin^{-1}(e^{\rhofin})\Jfin(e^{\rhofin+\alpha})
\end{equation*}
where $\Jfin = \sum_{\tau\in \Wfin} (-1)^\tau \tau$.
Let $q=e^{-\delta}$.  
Then for $\alpha\in\Pfin$,
\begin{equation} \label{shifted action}
 -\Pi\, e^{\ell\La_0} e^\alpha =
 \begin{cases}
   \Pi\,e^{\ell\La_0} e^{r_i(\alpha)-\alpha_i} & \text{for $i\not=0$} \\
   \Pi\, e^{\ell\La_0} e^{r_\theta(\alpha)+(\ell+1)\theta}
 q^{\ell+1-(\alpha|\theta)} & \text{for $i=0$.}
 \end{cases}
\end{equation}

These
equations express the $q$-equivalence in \cite{FLOT}.
Let $\Z^n$ have standard basis $\{\epsilon_i\mid 1\le i\le n\}$
and $\Pfin$ be the subspace of $\Z^n$ orthogonal to
the vector $\sum_{i=1}^n \epsilon_i$.  Then
$\alpha_i = \epsilon_i-\epsilon_{i+1}$ for $1\le i\le n-1$,
$\theta=\epsilon_1-\epsilon_n$, $(\cdot|\cdot)$ is the ordinary
dot product in $\Z^n$, and $\Wfin$ is the symmetric group on $n$
letters acting on the coordinates of $\Z^n$.
Since $\Pi \circ \Pifin = \Pi$ and $\Pifin$ is $\Z\La_0$-linear,
one may replace every term $e^\alpha$ by $s_\alpha := \Pifin e^\alpha$
in \eqref{shifted action}.
Define the map $\Z[\Pfin]^{\Wfin}[q]\rightarrow\Z[\Pfin]^{\Wfin}[q]$
given by $s_\alpha \mapsto \Pi(e^{\ell\La_0+\alpha}) e^{-\ell\La_0}$.
Define $f\equiv g$ in $\Z[\Pfin]^{\Wfin}[q]$ by the condition
that the above linear map sends $f$ and $g$ to the same element.
With this definition, we have
\begin{equation} \label{shifted equiv}
  -s_\alpha \equiv
  \begin{cases}
s_{(\alpha_1,\dots,\alpha_{i+1}-1,\alpha_i+1,\dots,\alpha_n)} &
\text{for $i\not=0$} \\
  s_{(\ell+1+\alpha_n,\alpha_2,\dots,\alpha_{n-1},-1-\ell+\alpha_1)}
  	q^{\ell+1-\alpha_1+\alpha_n} & \text{for $i=0$.}
\end{cases}
\end{equation}
It is not hard to see that this recovers the $q$-equivalence
of Schur functions given in \cite{FLOT}.

\subsection{Bosonic conjecture of \cite[(9.2)]{SW}}

In this section it is assumed that $\gggg$ is of type
$A^{(1)}_{n-1}$, $\La = \ell \La_0$, and the tensor factors $B_j$
are perfect crystals of the form $B^{k_j,\ell_j}$ in the
notation of \cite{KMN2} with $\ell_j\le \ell$ for all $j$.
By restriction to $U_q(\gfin)$, $B_j$ is the crystal of the
irreducible integrable $U_q(\gfin)$-module of highest weight
$\ell_j \Lab_{k_j}$.  In this case $B_0$ is not needed.
To see this, recall that $B_j$ can be realized as the
set of column-strict Young tableaux of the rectangular shape
having $k_j$ rows and $\ell_j$ columns with entries in the
set $\{1,2,\dots,n\}$.  In \cite{S2} the $\Pcl$-weighted
$I$-crystal structure on the perfect crystals
$B^{k,\ell}$ is computed explicitly.
In particular, if $b\in B_j$ is a tableau then
$\epsilon_0(b)$ is at most the number of ones in the tableau $b$,
which is at most $\ell_j$ by column-strictness.  Therefore
$b \otimes u_{\ell\La_0}$ never admits $e_0$.  Thus the
energy function $E_B$ of \eqref{energy} has the property that
for any $b\in B=B_L\otimes\cdots\otimes B_1$ such that
$e_0(b\otimes u_{\ell\La_0})=e_0(b)\otimes u_{\ell\La_0}$,
one has $E_B(e_0(b))=E_B(b)-1$.  Thus one obtains the
bosonic formula in this case.

Since $\gggg$ is of type $A^{(1)}_{n-1}$, $a_0=1$ and $\hck=n$. 
Take $\La=\La'=\ell\La_0$ in \eqref{level K poly}.
The lattice $M$ is given by
the root lattice $\Qfin$ of $\gfin$, which may be realized by
$\{\beta\in\Z^n\mid \sum_{i=1}^n \beta_i = 0\}$.
Let $B_{\tau,\beta}$ be the set of paths $b\in B$ of weight
$-\rhofin+\tau^{-1}(-(\ell+n)\beta+\rhofin)$.  Then
\begin{equation*}
\begin{split}
 \KP(B,\ell\La_0,\ell\La_0)(q) &=
 \sum_{\tau\in \Wfin} \sum_{\beta\in M} \sum_{b\in B_{\tau,\beta}}
 	(-1)^\tau q^{E_B(b)+(\rhofin\mid \beta)-
 	\tfrac{1}{2}\,|\beta|^2(\ell+n)} \\
 &= \sum_{\tau\in\Wfin}\sum_{\beta\in M} \sum_{b\in B_{\tau,\beta}}
 	(-1)^\tau q^{E_B(b)
         -\sum_{i=1}^n\{ \tfrac{1}{2}(\ell+n)\beta_i^2+i\beta_i\} }.
\end{split}
\end{equation*}
Notice that $\sum_{b\in B_{\tau,\beta}}q^{E_B(b)}$ is (up to an overall
factor) the $q\to1/q$ form of the supernomial $S$ of ref.~\cite{SW} so
that $K(B,\ell\La_0,\ell\La_0)(q)$ equals the left-hand side 
of~\cite[(9.2)]{SW} up to an overall power of $q$.
This shows that the left-hand side of~\cite[(9.2)]{SW} is indeed the 
generating function of level-$\ell$ restricted paths. To establish
the equality~\cite[(9.2)]{SW} it remains to prove that also the right-hand
side equals the generating function of level-restricted paths.

\subsection{Identities for level one and level zero}

As in the previous section let $\gggg$ be of type $A_{n-1}^{(1)}$
and assume that $B=B^{k_L,1}\otimes\cdots\otimes B^{k_1,1}$.
Fix $\ell=1$ and $\La,\La'\in\Pcl^+$ weights of 
level 1. It is easy to verify that $\Pu(B,\La,\La')$ consists of at most one 
element $p$. Choose $B,\La,\La'$ such that $p\in\Pu(B,\La,\La')$ exists.
Then by \eqref{K def} and \eqref{level K poly} we find that
\begin{equation} \label{id one}
 \sum_{\tau\in\Wfin}\sum_{\beta\in M}\sum_{b\in B_{\tau,\beta,\La,\La'}}
  (-1)^{\tau} q^{E(b)-\sum_{i=1}^n\{\tfrac{n+1}{2}\beta_i^2+i\beta_i\}}
 =q^{E(p)}
\end{equation}
where $B_{\tau,\beta,\La,\La'}$ is the set of paths $b\in B$ of
weight $-\Lab-\rhofin+\tau^{-1}(\Lab'-(n+1)\beta+\rhofin)$.

A similar formula exists for $\ell=0$:
\begin{equation} \label{id zero}
\sum_{\tau\in\Wfin}\sum_{\beta\in M}\sum_{b\in B_{\tau,\beta}}
  (-1)^{\tau} q^{E(b)-\sum_{i=1}^n\{\tfrac{n}{2}\beta_i^2+i\beta_i\}}
 =\delta_{B,\emptyset}
\end{equation}
where $B_{\tau,\beta}$ is the set of paths $b\in B$ of weight
$-\rhofin+\tau^{-1}(-n\beta+\rhofin)$. The right-hand side is the 
generating function of paths in $B$ of level zero since there are no
level zero restricted paths unless $B$ is empty. However, the arguments of 
Sections~\ref{sec general bose} and~\ref{sec inhomo paths}
do not imply that also the left-hand side is the generating function
of level zero paths since it was assumed in the proof of 
Theorem~\ref{gen boson} that the level of the crystals $B_j$ does not 
exceed $\ell$. We have assumed that $B_j=B^{k_j,1}$ which are crystals 
of level one. However, it is possible to define a sign-reversing involution
directly on $B=B^{k_L,1}\otimes\cdots\otimes B^{k_1,1}$ without using
the crystal isomorphisms that are used in the proof of Theorem~\ref{gen boson}.
Let $b\in B$. There exists at least one $0\le i\le n$ such that
$e_i(b_1)$ is defined. Define $v(b)=\min\{i|\text{$e_i(b_1)$ is defined}\}$
which has the property that $v(b)=v(\Phi_i(b))$ where as before 
$\Phi_i=s_ie_i$. Hence define the involution $\Phi(b)=\Phi_{v(b)}(b)$.
It is again sign-reversing and has no fixed points when $B\neq \emptyset$.
This proves that the left-hand side of~\eqref{id zero} is the generating
function of level 0 restricted paths.

Equation \eqref{id one} was conjectured in \cite{SW,W}.
For $n=2$ identity \eqref{id zero} follows from the $q$-binomial
theorem, for $n=3$ it was proven in~\cite[Proposition 5.1]{ASW}
and for general $n$ it was conjectured in~\cite{W}.

\end{document}